\documentclass[12pt]{article}
\newcommand{\authorfootnotes}{\renewcommand\thefootnote{\@fnsymbol\c@footnote}}%
\usepackage{amssymb,amsmath, amsfonts, amsthm}
\usepackage{verbatim,color}

\usepackage{graphicx}
\usepackage{tikz}
\usepackage{enumerate}
\usetikzlibrary{graphs}

\newcommand{\D}{\rm{Dominator\ }}
\newcommand{\St}{\rm{Staller\ }}

\newcommand{\cF}{\mathcal{F}}

\tikzstyle{vertex}=[circle, draw, inner sep=0pt, minimum size=6pt]

\usetikzlibrary{arrows}

\newtheorem{theorem}{Theorem}[section]
\newtheorem{lemma}[theorem]{Lemma}
\newtheorem{corollary}[theorem]{Corollary}
\newtheorem{proposition}[theorem]{Proposition}

\newtheorem{conjecture}[theorem]{Conjecture}

\newtheorem{obs}[theorem]{Observation}

\newcommand{\io}{\iota_{\rm g}}

\begin{document}

\title{Isolation game on graphs}

\author{
Boštjan Brešar$^{a,b}$
\and
Tanja Dravec$^{a,b}$
\and
Daniel P. Johnston$^{c}$
\and
Kirsti Kuenzel$^{c}$\\
\and
Douglas F.\ Rall$^{d}$
}

\date{\today}

\maketitle

\begin{center}
$^a$ Faculty of Natural Sciences and Mathematics, University of Maribor, Slovenia\\
$^b$ Institute of Mathematics, Physics and Mechanics, Ljubljana, Slovenia\\
$^c$ Department of Mathematics, Trinity College, Hartford, CT, USA\\
$^d$ Emeritus Professor of Mathematics, Furman University, Greenville, SC, USA\\
\end{center}

\begin{abstract}
Given a graph $G$ and a family of graphs $\cal F$, an $\cal F$-isolating set, as introduced by Caro and Hansberg, is any set $S\subset V(G)$ such that $G - N[S]$ contains no member of $\cal F$ as a subgraph. In this paper, we introduce a game in which two players with opposite goals are together building an $\cal F$-isolating set in $G$. Following the domination games, Dominator (Staller) wants that the resulting $\cal F$-isolating set obtained at the end of the game, is as small (as big) as possible, which leads to the graph invariant called the game $\cal F$-isolation number, denoted $\iota_{\rm g}(G,\cal F)$. We prove that the Continuation Principle holds in the $\cal F$-isolation game, and that the difference between the game $\cal F$-isolation numbers when either Dominator or Staller starts the game is at most $1$. Considering two arbitrary  families of graphs $\cal F$ and $\cal F'$, we find relations between them that ensure $\io(G,{\mathcal{F}}') \leq \io(G,{\mathcal{F}})$ for any graph $G$. A special focus is given on the isolation game, which takes place when $\cF=\{K_2\}$. We prove that $\io(G,\{K_2\})\le |V(G)|/2$ for any graph $G$, and conjecture that $\lceil 3|V(G)|/7\rceil$ is the actual (sharp) upper bound. We prove that  the isolation game on a forest when Dominator has the first move never lasts longer than the one in which Staller starts the game.
Finally, we prove good lower and upper bounds on the game isolation numbers of paths $P_n$, which lead to the exact values $\io(P_n,\{K_2\})=\left\lfloor\frac{2n+2}{5}\right\rfloor$
when $n \equiv i \pmod 5$ and $i \in \{1,2,3\}$.

\end{abstract}
\bigskip

\noindent
{\bf Keywords:} isolation number, graph game, domination game, Continuation Principle, forest \\

\noindent
{\bf AMS Subj.\ Class.\ (2020)}: 05C57, 05C05.

\maketitle

\section{Introduction}

Let $G=(V,E)$ be a graph and $\mathcal{F}$ a family of graphs. A graph $H$ is \emph{$\mathcal{F}$-forbidden} if no member of $\mathcal{F}$ is a subgraph of $H$. For a vertex $x\in V$, the {\em closed neighborhood} of $x$ is defined by $N_G[x]=\{v\in V:\, xv\in E\}\cup \{x\}$, and if $S\subset V$, the {\em closed neighborhood} of $S$ is $N_G[S]=\cup_{v\in S}{N[v]}$. A set $S$ of vertices is \emph{$\mathcal{F}$-isolating} if the subgraph $G-N_G[S]$ is $\mathcal{F}$-forbidden.  The minimum cardinality of an $\mathcal{F}$-isolating set in a graph $G$ is \emph{the $\mathcal{F}$-isolation number} of $G$ and is denoted by $\iota(G,{\mathcal{F}})$. The concept was introduced in 2017 by Caro and Hansberg~\cite{ch-2017}, and has gained the attention of many researchers. The $\mathcal{F}$-isolation number was studied for several families $\cal F$ and from different perspectives; see~\cite{bfk-2020,bk-2023,bg-2024,lms-2024,zw-2024} for a short selection of such studies.

In this paper, we introduce a game counterpart of the isolation number. The idea for the game was inspired by the domination game, introduced in 2010 by Bre\v{s}ar, Klav\v{z}ar and Rall~\cite{bkr-2010} as a game counterpart of the domination number. Numerous papers were devoted to the study of domination games, and several innovative proof methods were introduced. In particular, the Imagination Strategy and the Continuation Principle are two such important tools; see~\cite{bhkr-2021} and the references therein. In addition, a special discharging method due to B\'{u}jtas~\cite{bujtas} was developed while investigating the domination game, and was later applied for obtaining best upper bounds for graphs with specified minimum degree on the domination number itself~\cite{BK-2016}. Many other variations of the domination game were introduced, each presenting a game version of a graph invariant that may be recognized by its name: total domination game~\cite{hkr-2015}, fractional domination game~\cite{bt-2019}, connected domination game~\cite{bfs-2019}, two variants of the paired-domination game~\cite{gh-2023,hh-2019}, and competition-independence game~\cite{ph-sl-1,ph-sl-2}.
More examples of such variations can be found in the book~\cite{bhkr-2021} surveying the domination game and related concepts.

The \emph{$\mathcal{F}$-isolation game} is played on a graph $G$ by two players, Dominator and Staller, who take turns choosing (playing) a vertex from $G$  while obeying the following rule. If $S$ is the set of already played vertices, then vertex $x$ can be chosen only if it dominates a vertex $y$ that belongs to a component of the graph $G-N[S]$ that is not $\mathcal{F}$-forbidden. (Recall that a vertex $x$ {\em dominates} itself and all its neighbors.) Such a vertex $x$ is called {\em playable}, and a move selecting $x$ is {\em legal}. The game ends when no playable vertex exists. Note that when the game has ended, the set of selected vertices is an $\mathcal{F}$-isolating set. The players have opposite goals.
Dominator wishes to finish the game as fast as possible, and Staller wishes to delay the process as much as possible.
The {\em game $\mathcal{F}$-isolation number}, denoted by $\io(G,{\mathcal{F}})$, is the number of vertices chosen when Dominator starts the game
(the {\em Dominator-start game}) provided that both players play optimally. Similarly,  the {\em Staller-start game $\mathcal{F}$-isolation number}, written as $\io'(G,\mathcal{F})$, is the number of vertices chosen in
 the game when Staller starts the game (the {\em Staller-start game}). The Dominator-start game and the Staller-start game will be briefly
called the {\em D-game} and the {\em S-game}, respectively, when the family $\mathcal{F}$ is clear from the context.

If $\mathcal{F}=\{K_1\}$, then the $\mathcal{F}$-isolation game coincides with the well-known domination game. Thus $\io(G,\{K_1\})=\gamma_{\rm{g}}(G)$ and $\io'(G,\{K_1\})=\gamma'_{\rm{g}}(G)$.
If $\mathcal{F}=\{K_2\}$, then the vertices not dominated by the $\mathcal{F}$-isolating set form an independent set. In this case the $\mathcal{F}$-isolation number is denoted by $\iota(G)$, and we will use the notation $\io(G)$ instead of $\io(G,{\mathcal{F}})$ and $\io'(G)$ instead of $\io'(G,{\mathcal{F}})$.

In the next section, we consider the $\cal F$-isolation game from a general perspective, providing results in which the family $\cal F$ is arbitrary. By using the Imagination Strategy we prove that the Continuation Principle holds for all $\cF$-isolation games regardless of the choice of $\cal F$. Using this tool, we then show that $|\io(G,{\mathcal{F}})-\io'(G,{\mathcal{F}})|\leq 1$, where $G$ is an arbitrary graph and $\cal F$ an arbitrary family of graphs. Then, we prove a simple lower and an upper bound on the game $\cal F$-isolation number of a graph $G$ in terms of the $\cF$-isolation number of $G$, and provide a certain ``hereditary'' condition for two families $\cF$ and $\cF'$, which ensures that $\io(G,{\mathcal{F}}') \leq \io(G,{\mathcal{F}})$. In Section~\ref{sec:iso}, we focus on the isolation game (that is, when $\cF=\{K_2\}$), and prove that $\io(G)$ is bounded from above by half the order of the graph $G$. We conjecture that the bound can be improved to $\io(G)\le\lceil\frac{3n(G)}{7}\rceil$, where $n(G)$ is the order of $G$, and present an infinite family of graphs that attain this bound. Section~\ref{sec:forest} is devoted to the isolation game on forests. We prove that the D-game in a forest never lasts longer than the S-game (if both players play optimally). We also provide a lower and an upper bound on the game isolation number of a path, where the bounds are at most $1$ apart.  When $n \equiv i \pmod 5$, where $i \in \{1,2,3\}$, we get the exact values: $\io(P_n)=\io'(P_n)= \left\lfloor\frac{2n+2}{5}\right\rfloor$.

\section{$\cal F$-isolation game}

In this section, we collect a number of preliminary results, most of which can be considered as extensions from the domination game to the isolation game. The proofs of several of these results follow similar lines of reasoning as the proofs of corresponding results for the domination game. For the sake of completeness we will provide all these proofs.

During the course of the game, when $S$ is the set of selected vertices, a vertex $v$ is considered \emph{marked} if $v \in N[S]$ or $v$ is contained in an $\mathcal{F}$-forbidden component of $G - N[S]$. In particular, before the game begins every vertex in an $\cF$-forbidden component, if such a component exists, is marked. Note that a marked vertex may or may not be playable. It is easy to see the following observation about playable vertices during the course of the game.
\begin{obs}
\label{obs:playable}
 When the $\cF$-isolating game is played on a graph, a vertex is not playable if and only if its entire closed neighborhood is marked.
\end{obs}

A \emph{partially marked graph} $G|A$ consists of  a graph $G$ with a subset $A\subseteq V(G)$ such that $G-A$ does not contain an $\mathcal{F}$-forbidden component. We will consider the vertices in $A$ as already marked. Suppose we are playing the game on a partially marked graph $G|A$. During the course of the game, after each move we will update the set of marked vertices $M$ where initially $M=A$. Suppose $x$ is a playable vertex in $G|M$. Let \[U = \{t: t \text{ is a vertex in an $\mathcal{F}$-forbidden component of $G - (M \cup N[x])$}\},\]
and let $M' = M \cup U \cup N[x]$.
If a player selects the vertex $x$ when the game is played on $G|M$, then $G|M'$ is the updated partially marked graph. The game ends when every vertex in
$G$ is marked.

One of the main tools in analyzing the domination game is the Continuation Principle~\cite{KWZ-2013}. In the next result we show that this principle extends to the $\mathcal{F}$-isolation game.

\begin{lemma}\label{l:continuation}{\em({Continuation Principle)}}
    If $G$ is a graph and $A,B \subseteq V(G)$ with $B \subseteq A$, then $\io(G|A,{\mathcal{F}}) \leq \io(G|B,{\mathcal{F}})$ and $\io'(G|A,{\mathcal{F}}) \leq \io'(G|B,{\mathcal{F}})$, where $\cF$ is an arbitrary family of graphs.
\end{lemma}

\proof

Two games will be played in parallel, Game~A on the partially marked graph $G|A$ and Game~B on the partially marked graph $G|B$. The first of these will be the real game, while Game~B will be imagined by Dominator. In Game~A, Staller will play optimally while in Game~B, Dominator will play optimally.

 We claim that in each stage of the games, the set of vertices that are marked in Game~B is a subset of the set of vertices that are marked in Game~A. Since $B \subseteq A$, this is true at the start of the games. Suppose now that Staller has (optimally) selected a vertex $u$ in Game~A. Then by the induction assumption, the vertex $u$ is playable in Game~B because a new vertex $v$ that was dominated by $u$ in Game~A, is not yet marked in Game~B as well. Then Dominator copies the move of Staller by playing the vertex $u$ in Game~B. \D then replies with an optimal move selecting vertex $x$  in Game~B. (In particular, if the D-game is being played and this is the first move in the game, then \D will choose an optimal move in Game~B.) If $x$ is playable in Game~A, \D plays $x$ in Game~A as well. Otherwise, if the game is not yet over, \D selects any other playable vertex in Game~A. In either case the set of marked vertices  in Game~B is a subset of the set of marked vertices in Game~A which by induction also proves the claim.

We have thus proved that Game~A finishes no later than Game~B. Suppose that Game~B lasted $r$ moves. Because Dominator was playing optimally in Game~B, it follows that $r \le \io(G|B, \mathcal{F})$. On the other hand, because Staller was playing optimally in Game~A and \D has a strategy to finish the game in $r$ moves, we infer that $\io(G|A, \mathcal{F})\le r$. Therefore, $\io(G|A, \mathcal{F})\le r \le  \io(G|B, \mathcal{F})$ and we are done if \D is the first to play. Note that in the above argument we did not assume who starts first, hence in both cases Game~A will finish no later than Game~B. Hence the conclusion holds for $\io'$ as well; that is, $\io'(G|A, \mathcal{F})\le \io'(G|B, \mathcal{F})$.
\qed

\medskip

As a consequence of the Continuation Principle, whenever $x$ and $y$ are playable vertices for \D in the $\mathcal{F}$-isolation game and $N[x] \subseteq N[y]$, we may assume that Dominator will play $y$ instead of $x$, since any newly marked vertex as a result of playing $x$ will also be marked as a result of playing $y$.  For a similar reason, if it was Staller's turn, we may assume that Staller will play $x$ instead of $y$.

\begin{theorem}\label{thm:difference1}
    For an arbitrary graph $G$ and an arbitrary family $\cF$, we have $|\io(G,{\mathcal{F}})-\io'(G,{\mathcal{F}})|\leq 1$.
\end{theorem}

\proof
Consider the D-game on $G$ and let $v$ be the first move of Dominator.
Let $A$ be the set of marked vertices of $G$ after $v$ is played, and let $B=\emptyset$. Consider the partially marked graph $G|A$, and note that $G|B = G$. Since $v$ is an optimal first move of Dominator, we note further that $\io(G, \mathcal{F}) = 1 + \io'(G|A, \mathcal{F})$. By the Continuation Principle, $\io'(G|A, \mathcal{F}) \le \io'(G|B, \mathcal{F}) = \io'(G, \mathcal{F})$.
Therefore,
\[
\io(G, \mathcal{F}) = \io'(G|A, \mathcal{F}) + 1 \le \io'(G, \mathcal{F}) + 1\,.
\]
By a parallel argument, consider the S-game and let $v$ be the first move of Staller. As before, let $A$ be the set of marked vertices of $G$ after $v$ is played by \St and $B=\emptyset$. Since $v$ is an optimal first move of Staller, we note that $\io'(G, \mathcal{F}) = 1 + \io(G|A, \mathcal{F})$. By the Continuation Principle, $\io(G|A, \mathcal{F}) \le \io(G|B, \mathcal{F}) = \io(G, \mathcal{F})$, implying that
\[
\io'(G, \mathcal{F}) = \io(G|A, \mathcal{F}) + 1 \le \io(G, \mathcal{F}) + 1\,.
\]
\vskip -0.6cm
.~\qed

\medskip

\begin{theorem}
    For an arbitrary graph $G$ and an arbitrary family $\cF$, we have $$\iota(G,{\mathcal{F}}) \leq \io(G,{\mathcal{F}}) \leq 2\iota(G,{\mathcal{F}}) - 1\textrm{,  and } \iota(G,{\mathcal{F}}) \leq \iota'_g(G,{\mathcal{F}}) \leq 2\iota(G,{\mathcal{F}}).$$
\end{theorem}
\proof
Since for any graph $G$ and any family $\cF$ the set of selected vertices in the $\cF$-isolation game played on $G$ is an $\cF$-isolating set, we get $\iota(G,{\mathcal{F}}) \leq \io(G,{\mathcal{F}})$. Let $I$ be an $\cal F$-isolating set of $G$, where $\iota(G,\cF)=|I|$. Dominator's strategy when the D-game is played on $G$ is to select the vertices from $I$ one at a time as long as any of these vertices is playable.  When no additional vertex from $I$ is playable, we infer by Observation~\ref{obs:playable} that every vertex in $N[I]$ is marked, which by the choice of $I$ means that all vertices of $G$ are marked. Thus, altogether at most $2|I|-1$ moves have been played in the game, yielding $\io(G,{\mathcal{F}}) \leq 2\iota(G,{\mathcal{F}}) - 1$.

The inequalities for $\io'(G,\cF)$ can be proved by using similar reasoning.
\qed

\medskip

The next result relates isolating games on the same graph $G$ considering two different families $\cal F$ and $\cF'$, where $\cal F$ and $\cF'$ are related in a specific way.

\begin{theorem}\label{thm:relations}
Let $\mathcal{F}$ and $\mathcal{F}'$ be two families of graphs. If for any $F \in {\mathcal{F}}$ there exists $F' \in {\mathcal{F}}'$ such that $F$ is a subgraph of $F'$, then $\io(G,{\mathcal{F}}') \leq \io(G,{\mathcal{F}})$, where $G$ is an arbitrary graph.
\end{theorem}
\proof
To prove $\io(G,{\mathcal{F}}') \leq \io(G,{\mathcal{F}})$, we will consider the $\cF'$-isolation D-game played on $G$, while Dominator will imagine that an $\cF$-isolation D-game is played in parallel on $G$. We call these games the $\cF'$-game and the $\cF$-game, respectively.  Let $A'$, respectively $A$, be the set of marked vertices updated after every move of the players, in the $\cF'$-game, respectively $\cF$-game on $G$. Note that, by the assumption on the families $\cF$ and $\cF'$, initially $A\subseteq A'$.
By the Dominator's strategy that we will present next, the relationship between $A$ and $A'$ will be maintained after each move.

When Dominator selects $v$ in the $\cF$-game, then he plays $v$ also in the $\cF'$-game if $v$ is playable there. Otherwise, he chooses any playable vertex in the $\cF'$-game.
Clearly, in either case, after his move the relation $A\subseteq A'$ still holds. Staller is following her optimal strategy in the $\cF'$-game. Each of her moves in this game is copied by Dominator to the $\cF$-game, and since $A\subseteq A'$, this move is also legal in the $\cF$-game. Because of the assumed relationship between the families $\cF$ and $\cF'$, the updated set of marked vertices in the $\cF$-game is a subset of the updated vertices in the $\cF'$-game. Let $s$ and $s'$ be the total numbers of moves in the $\cF$-game and the $\cF'$-game, respectively. By the relationship between the marked vertices in the two games, we infer $s'\le s$. Since Staller played optimally in the $\cF'$-game, it follows that $\io(G,\cF')\le s'$. Similarly, since Dominator  followed his optimal strategy in the $\cF$-game, we infer $s\le \io(G,\cF)$.
\qed

\medskip

Since $\io(G)=\io(G,\{K_2\})$ and $\gamma_{\rm g}(G)=\io(G,\{K_1\})$ Theorem~\ref{thm:relations} implies the following.

\begin{corollary}\label{cor:relationToGamma}
    For any graph $G$ we have $\io(G) \leq \gamma_{\rm g}(G)$.
\end{corollary}

Note that the difference $\gamma_{\rm g}(G)-\io(G)$ can be arbitrarily large.
Let $n$ be any positive integer and let $T_1, \ldots, T_n$ be $n$ vertex disjoint triangles, where $V(T_i)=\{v_i,x_i,y_i\}$ for each $i \in [n]$.  Let $G_n$ be the graph constructed
from the union of this set of triangles by letting $\{v_1,\ldots,v_n\}$ induce a complete subgraph.  In addition, for each $k \in [n]$ let
$F_k$ be the spanning subgraph of $G_n$ defined by $F_k=G_n -\{x_1y_1, \ldots, x_ky_k\}$.
Now, letting $n$ be odd, it is easy to show that $\io(F_n)=1$ and $\gamma_{\rm g}(F_n)=1+3(n-1)/2$.

\section{Isolation game}
\label{sec:iso}

In the rest of the paper $\cF=\{K_2\}$, that is, we consider the isolation game.  As mentioned earlier, the game $\{K_2\}$-isolation invariants are denoted by $\io$ and $\io'$.

In the next example, we show that the game isolation number for a spanning subgraph $H$ of $G$ can be smaller than $\io(G)$. Let $G_n$ and its spanning subgraphs $F_k$, $k \leq n$, be the graphs as defined in the previous section. Then $\io(G_n)=n$ and
$\io(F_k)=n-k$, for each $k \in [n-1]$ and $\io(F_n)=1$.  This example shows that not only can the game isolation number of a spanning subgraph be arbitrarily smaller than the game isolation number of the original graph, but it also shows that the game isolation number can ``interpolate'' between the smallest and largest values. To see that $\io(H)$ can also be arbitrarily larger than $\io(G)$, let $G=K_n$. We summarize this discussion as follows.

\begin{obs}
\label{obs:spanning}
For an arbitrary $k\in \mathbb{Z}$, there exists a graph $G$ with a spanning subgraph $H$ such that $\io(G)-\io(H)=k$.
\end{obs}

In the next result, we provide a general upper bound on the game isolation number of a graph in terms of a fraction of its order.

\begin{theorem}\label{thm:1/2bound}
    For any graph $G$, $\io(G) \leq \frac{1}{2}n(G)$.
\end{theorem}

\proof
We will weight the graph in the following way. At the beginning all vertices are unmarked and each vertex has weight 1. During the game, when a vertex is marked its weight drops to 0. The weight of $G$, $w(G)$, is the sum of weights of all vertices of $G$. Thus at the beginning of the game $w(G)=n(G)$ and when the game ends the weight is zero. Let $M$ be the set of marked vertices in some stage of the game. Dominator's strategy is that he does not play a vertex in any $K_2$ component of $G-M$, unless he has no other option. Moreover, if there exists a component $C$ of $G-M$ of order at least 3, then he plays a vertex of degree at least 2 in $C$. With this strategy,  each of his moves  will decrease the weight of $G$ by at least 3. Since Staller will decrease the weight of $G$ by at least 1 after each of her moves, each move will on average decrease the weight by at least 2 until just $K_2$ components are left in $G-M$. In this final stage of the game, if any $K_2$ components are left in $G-M$, each player will in every move decrease the weight by exactly 2. Hence the game finishes after at most $\frac{n(G)}{2}$ moves.

\qed

The bound from Theorem~\ref{thm:1/2bound} is sharp since $\io(C_6)=3$. However, we think that the bound can be improved to $\io(G) \leq \left\lceil\frac{3}{7}n(G)\right\rceil$. We could not find an example where this bound is exceeded, and we pose the following conjecture.

\begin{conjecture}\label{con:3/7}
    For any graph $G$ with no $K_2$-components, $$\io(G) \leq \left\lceil\frac{3}{7}n(G)\right\rceil \text{ and } \io'(G) \leq \left\lceil\frac{3}{7}n(G)\right\rceil.$$
\end{conjecture}

The ceiling in the expression $3n(G)/7$ is needed, since $\io(C_6)=3$ and $\io'(P_4)=2$.

Next, we present an infinite family of graphs showing that the conjectured bound can be attained. Let $G$ be an arbitrary graph. For each vertex $v \in V(G)$ identify $v$ with the center vertex of a path of order $7$ and denote the resulting graph by $G^*$.

\begin{proposition}
\label{prp:familyequality}
  If $G$ is an arbitrary graph, then $\io(G^*)=\io'(G^*)= \frac{3}{7}n(G^*)$.
\end{proposition}
\proof
    Let us call the paths of order $7$ that are in $G^*$ but not in $G$ ``added paths.''  Let D-game be played on $G$. A strategy of Dominator that no more than three vertices are played on any of the added paths is as follows. If he is the first player to select a vertex on an added path, then he achieves this by choosing the center of the path. (Indeed, after his move in the added path in which he played, two $K_2$ components of unmarked vertices remain, and there are exactly two moves left in this added path.) Now suppose that Staller is the first player to select a vertex from an added path.  If she does not choose the center vertex of that path, then Dominator can prevent four played vertices from that path by choosing a vertex that is at distance 4 from her move. On the other hand, Staller has a strategy to force three moves on each added path. If she has the first move on some added path, then she chooses the center; otherwise she responds to Dominator's first move on this path by choosing its neighbor.

Following the same strategies also in the S-game, one can easily derive $\io'(G^*)=\frac{3}{7}n(G^*)$.
\qed
\bigskip

Therefore, if Conjecture~\ref{con:3/7} holds, then it is sharp for an infinite class of graphs.
Besides the family of graphs $G^*$ that appear in Proposition~\ref{prp:familyequality} and graphs similar to graphs $G^*$ (say, obtained from $G^*$ by removing a pendant vertex), we could not find any other infinite family of graphs whose game isolation number is close to $3n/7$. For instance, for paths $P_n$ (see Section~\ref{sec:forest} for details) we prove that the game isolation number is roughly $2n/5$. After searching for graphs having game isolation number as large  as possible, we eventually found an infinite family $\cal{G}$ of graphs with $\frac{2n(G)}{5} < \io(G)=\frac{5n(G)}{12}  < \frac{3n(G)}{7}$ for all $G \in \cal{G}$, which we present next.

\begin{figure}[ht!]
\begin{center}
\begin{tikzpicture}[scale=0.6,style=thick,x=0.9cm,y=0.9cm]
\def\vr{5pt}
\path (4,1) coordinate (v1); \path (6,2) coordinate (v2); \path (6,0) coordinate (v3); \path (2,1) coordinate (v4);
\path (0,2) coordinate (v5); \path (0,0) coordinate (v6); \path (0,-2) coordinate (v7); \path (-1,-4) coordinate (v8);
\path (1,-4) coordinate (v9); \path (0,4) coordinate (v10); \path (-1,6) coordinate (v11); \path (1,6) coordinate (v12);

\draw (v1) -- (v2); \draw (v1) -- (v3);  \draw (v1) -- (v4);  \draw (v2) -- (v3);  \draw (v4) -- (v5);  \draw (v4) -- (v6);
\draw (v5) -- (v6); \draw (v5) -- (v10);  \draw (v6) -- (v7);  \draw (v7) -- (v8);  \draw (v7) -- (v9);  \draw (v8) -- (v9);
\draw (v10) -- (v11); \draw (v10) -- (v12);  \draw (v11) -- (v12);

\foreach \i in {1,...,12}
{  \draw (v\i)  [fill=white] circle (\vr);}
\draw (v4) [fill=black] circle (\vr);

\draw (4,1.6) node {{$v_1$}}; \draw (6,2.6) node {{$v_2$}}; \draw (6,-.6) node {{$v_3$}}; \draw (2,1.6) node {{$v_4$}}; \draw (-.6,2) node {{$v_5$}};
\draw (-.6,0) node {{$v_6$}}; \draw (-.6,-2) node {{$v_7$}}; \draw (-1.6,-4) node {{$v_8$}}; \draw (1.6,-4) node {{$v_9$}}; \draw (-.6,4) node {{$v_{10}$}};
\draw (-1,6.6) node {{$v_{11}$}};\draw (1,6.6) node {{$v_{12}$}};

\end{tikzpicture}
\caption{Graph $H$ with $\io(H)=5=\io'(H)$}
\label{fig:example}
\end{center}
\end{figure}

For a positive integer $n$ let $G$ be an arbitrary graph of order $n$ with $V(G)=\{w_1,\ldots,w_n\}$ along with $n$ vertex disjoint copies,
say $H_1, \ldots, H_n$, of the graph $H$ depicted in Figure~\ref{fig:example}.  For each $i \in [n]$, identify the vertex $w_i$ in $G$ with the vertex $v_4$ in $H_i$, and call the resulting graph $G_H$.

One can check that for the graph $H$ in Figure~\ref{fig:example} it holds $\io(H|\{v_4\})=5=\io'(H|\{v_4\}) \text{ and } \io(H)=5=\io'(H)$. Furthermore, it is not hard to verify that skipping a move when playing on $H$ does not help either of the players. Using these facts, one can then find strategies for both players in the D-game and in the S-game on $G_H$ that lead to exactly five moves being played on each $H_i$. Thus, for every positive integer $n$, we have $\io(G_H)=\io'(G_H)=5n = \frac{5|V(G_H)|}{12}$, where $G$ is an arbitrary graph of order $n$.

\section{Isolation game on forests}
\label{sec:forest}

The proofs of the  next lemma and theorem closely follow arguments of Kinnersley, West, and Zamani~\cite{KWZ-2013} for the domination game.  They are included for completeness.
We write $G + H$ to denote the disjoint union of graphs $G$ and $H$. If $G$ is a partially marked graph and $v \in V(G)$, then $G_v$ denotes the partially marked graph obtained from $G$ by marking all vertices corresponding to a play of $v$.
\begin{lemma}\label{l:forests}
    Given $\lambda \geq 2$, assume $\io(F) \leq \io'(F)$ for all partially marked forests $F$ such that $\io(F) \leq \lambda$.  If $G$ is such a forest and $r$ is any non-negative integer, then $\io(G + K_{1,r}) > \io(G)$ and $\io'(G + K_{1,r}) > \io'(G)$.
\end{lemma}

\proof
Fixing $\lambda$, we use induction on the number of  unmarked vertices of $G$.  The claim is true if $G$ has no unmarked vertices since an additional play will be needed to mark the vertices of $K_{1,r}$.

To prove $\io(G + K_{1,r}) > \io(G)$, let $v$ be an optimal first move in the D-game on $G + K_{1,r}$.  If $v$ is a vertex of $K_{1,r}$, then $\io(G + K_{1,r}) = 1 + \io'(G) \geq 1 + \io(G)$.  Otherwise, $v$ belongs to $G$.  By the optimality of $v$ we have $\io(G+K_{1,r}) = 1 + \io'(G_v+K_{1,r})$, and by the Continuation Principle, $\io(G_v) \leq \io(G) \leq \lambda$.  These observations and the induction hypothesis yield $\io(G+K_{1,r}) = 1 + \io'(G_v+K_{1,r}) > 1 + \io'(G_v)$.  Since $v$ is not necessarily an optimal first move for Dominator in the D-game on $G$, we have $1+ \io'(G_v) \geq \io(G)$ and hence $\io(G+K_{1,r}) > \io(G)$.

To prove $\io'(G + K_{1,r}) > \io'(G)$, let $v$ be an optimal first move in the S-game on $G$.  Now $\io'(G + K_{1,r}) \geq 1 + \io(G_v + K_{1,r}) > 1 + \io(G_v) = \io'(G)$ with the strict inequality following from the induction hypothesis and the equality following from the optimality of $v$.
\qed

\medskip

As for the game domination number, there exist graphs $G$ where $\io(G) < \io'(G)$, graphs $G$ with $\io(G)=\io'(G)$ (one such example is $P_3$) and graphs $G$ with $\io(G) > \io'(G)$ (one such example is $C_6$, for which $\io(C_6)=3$ and $\io'(C_6)=2$). In the next theorem we will prove that only two of these three relationships are  possible when $G$ is a forest.

\begin{theorem}\label{thm:MonotonicityForests}
    If $F$ is a partially marked forest, then $\io(F) \leq \io'(F)$.
\end{theorem}
\proof
The proof goes by induction on $\io(F)=k$. If $k\leq 2$, the result trivially follows.  Indeed, if $\io(F)=2$, then it is not possible that $\io'(F) \leq 1$. Thus let $k \geq 3$ and assume that there exists a partially marked forest $F$ with $\io(F)=k$ and $\io'(F) \leq k-1$. By Theorem~\ref{thm:difference1}, $\io'(F)=k-1$. Therefore,  Dominator can ensure that after any first move of Staller in the S-game on $F$, the game ends within at most $k-2$  additional moves. That is, $\io(F_v) \leq k-2$ for any $v \in V(F)$. Similarly, since $\io(F)=k$, we get $\io'(F_v) \geq k-1$ for any $v \in V(F)$. By Theorem~\ref{thm:difference1}, $\io(F_v) = k-2$ and $\io'(F_v) =k-1$ hold for any $v \in V(F)$.

Denote by $M$ the set of vertices of $F$ that are marked before the game starts. Note that by definition each vertex that is not marked has at least one unmarked neighbor. Let $C$ be a component of $F$  having an unmarked vertex. We root $C$ in an arbitrary vertex $r$ and let $x$ be an unmarked vertex in $C$ that is at the largest distance from $r$. Since $x \notin M$ and all descendants of $x$ are in $M$, the only possible unmarked neighbor of $x$ is its parent $x'$. If $r\neq x'$, let $y$ be the parent of $x'$, otherwise let $y=x'$. Let $u$ be an optimal first move of Staller in the S-game played on $F_y$. By the choice of $x$ and by the definition of $F_y$, it follows that $u \notin N[\{x,x'\}]$. Since $u$ is an optimal first move of Staller in the S-game on $F_y$, $\io'(F_y)=k-1$ implies that $\io((F_y)_u)=k-2$. Moreover, we already know that $\io(F_u)=k-2$.

Let $F'$ be the partially marked forest obtained from $F$ by marking $X=N[u]\cup (N[\{y,x'\}]\setminus \{x,x'\})$ and all isolated vertices of the graph $F-(M \cup X)$. Note that the only difference between $F'$ and $(F_y)_u$ is  that vertices $x,x'$ are marked in $(F_y)_u$  but $x$ and $x'$ are not marked in $F'$. Since the move on $F'$ that marks either $x$ or $x'$ marks both $x$ and $x'$ but no other vertices of $F'$, we get $\io(F')=\io((F_y)_u+K_{1,1})$, where $K_{1,1}$ corresponds to the star on vertices $x'$ and $x$. Hence $$ k-2 = \io(F_u) \geq \io(F')=\io((F_y)_u+K_{1,1}) > \io((F_y)_u)=k-2,$$
which is a contradiction. Note that the first inequality follows from the Continuation Principle and the last inequality follows from Lemma~\ref{l:forests}.
\qed

\medskip

By Theorem~\ref{thm:MonotonicityForests},  we can obtain a good upper bound for the game isolation number of paths.

\begin{lemma} \label{lemma:DRupperbound}
For any $n \geq 3$,
    $\io(P_n) \le \io'(P_n)\le \left\lfloor\frac{2n+2}{5}\right\rfloor$
\end{lemma}
\proof
The first inequality follows from Theorem~\ref{thm:MonotonicityForests}. For the second inequality consider the S-game played on $P_n$.  Dominator can use the strategy following each of Staller's moves by selecting a vertex $x$ that is at distance $4$ from one of the already played vertices, say $y$, provided there are no vertices lying between $x$ and $y$ that have been played. In this way, if such a vertex $x$ is played, Dominator's move adds  four vertices to the set of unplayable vertices. Clearly, each move of Staller makes at least one vertex unplayable. Suppose that Dominator can no longer select such a vertex as described above. Then, note that there are exactly three options:
\begin{itemize}
    \item the game is over (all vertices are unplayable), or
    \item there are exactly three playable vertices, which lie in one end of the path, or
    \item there are exactly six playable vertices, three at each end of the path.
\end{itemize}

Let $u_i$ be the number of vertices that became unplayable because of the $i$th move in the game, and let $p$ be the total number of moves made in the game. In the first case, and we get
\[
n = \sum_{i=1}^p u_i \geq
\left\{\begin{array}{lr}
5\frac{p}{2};  \text{ Dominator had the last move},\\
5\frac{p-1}{2}+3;  \text{ Staller had the last move},
\end{array}\right.
\]
both of which imply that $ p \le \frac{2n}{5}$.  Since $p$ is an integer we get $p \leq \lfloor \frac{2n}{5}\rfloor$. In the second case, the last move, say $w$, that resulted in exactly three playable vertices
was played by Staller due to Dominator's strategy. Hence
\[
n = \sum_{i=1}^p u_i \geq 5\frac{p-2}{2} +1+3\,
\]
and this gives  $p \leq \frac{2n+2}{5}$.  Again, since $p$ is an integer, it follows that $p \leq \lfloor \frac{2n+2}{5}\rfloor$.

Finally, in the last case the last move, say $w$, that resulted in exactly six playable vertices was again played by Staller for the same reason as in the previous case.   Therefore,
\[
n = \sum_{i=1}^p u_i \geq 5\frac{p-3}{2} +1+6\,
\]
and we get $p \leq \frac{2n+1}{5}$ and so $p \leq \lfloor \frac{2n+1}{5}\rfloor$.

We have shown that by following the strategy described above, Dominator can ensure that $p \leq \lfloor \frac{2n+2}{5}\rfloor$.
Therefore, $\io'(P_n)\le p \leq \left\lfloor\frac{2n+2}{5}\right\rfloor$.

\qed

\bigskip

We also get a good lower bound:

\begin{lemma}\label{l:lowerBound}
For any $n \geq 6$,
    $ \lceil\frac{2n}{5}\rceil -1\le \io(P_n) \le \io'(P_n).$
\end{lemma}
\proof
We only need to prove $\lceil\frac{2n}{5}\rceil -1\le \io(P_n)$.
Note that when the D-game played on a path is over, the following sets $A_1,\ldots,A_{k+1}$, and $B_1,\ldots, B_k$ are generated. If $v_1$ was selected during the game, then $A_1=\emptyset$, and $B_1$ is a maximal string of vertices that were selected during the game, starting with $v_1$. Otherwise, $A_1$ is a maximal string of vertices that were not selected during the game including $v_1$, which ends with, say $v_{k_1}$. In that case, $B_1$ is a maximal string of vertices that were selected during the game, starting with $v_{k_1+1}$. We continue in this way, where finally $A_{k+1}$ may or may not be empty. Observe that $|A_1|\le 2$, $|A_{k+1}|\le 2$, and $|A_i|\le 3$ for $i\in\{2,\ldots,k\}$ regardless of how the game is played.

Staller's strategy to achieve the desired bound is as follows. If Dominator selected a vertex $v$ in one of his moves, then Staller responds by selecting a neighbor of $v$ if possible. (This is clearly possible for her after the first move of Dominator.)  If this is not possible, then Staller selects any vertex that is a neighbor of a vertex selected earlier. If there are no such vertices available, the game is already over. Note that after each move of Dominator (with the only possible exception if Dominator played the last move in the game), Staller's strategy of selecting a vertex $u$ next to another already selected vertex $v$ increases the size of the corresponding set $B_\ell$, in which the vertex $u$ lies, by $1$. In particular, if $v$ was played in the preceding move by Dominator, then after the move of Staller selecting $u$ we have $|B_\ell|\ge 2$, while if $v$ was played in an earlier move, then Staller's selection of $u$ yields $|B_\ell|\ge 3$. Based on these facts, we have the following correspondence at the end of the game: for every set $B_i$ with $|B_i|=1$ (with only one possible exception when the vertex $v$ of $B_i$ was the last move in the game played by Dominator), where $v$ is the vertex of $B_i$ necessarily played by Dominator,  there exists a set $B_j$ with $|B_j|\ge 3$ and a vertex $u$ of $B_j$, which was played by Staller as the response to Dominator's selection of vertex $v$.

Let $p$ be the number of moves played in the D-game when Staller employs the above strategy. We infer $p= \sum_{i=1}^k{|B_i|}$, which is either $2k$ or $2k-1$ depending on whether Dominator had the last move in the game and with this move selected a vertex $v$, which belongs to $B_i$ with $|B_i|=1$. On the other hand, $$\sum_{i=1}^{k+1}{|A_i|}\le 2+2+3(k-1)=3k+1.$$
First, if $p=2k$, then $n\le 5k+1$, which gives $p\ge \frac{2n-2}{5}$. On the other hand, if $p=2k-1$, then $n\le 5k$, which gives $p\ge \frac{2n-5}{5}$. In either case, we get that the number of moves is at least $\frac{2n}{5}-1$, and since this is an integer,
$\lceil\frac{2n}{5}\rceil -1\le \io(P_n)$, as claimed.
\qed

\bigskip

Combining the upper bound from Lemma~\ref{lemma:DRupperbound} and the lower bound from Lemma~\ref{l:lowerBound} we get the following.

\begin{corollary}\label{cor:DRpaths}
If $n \geq 6$,  then
 \[\left\lceil\frac{2n}{5}\right\rceil-1\le \io(P_n)\le\io'(P_n) \le  \left\lfloor\frac{2n+2}{5}\right\rfloor.
 \]
\end{corollary}

\medskip

Using Corollary~\ref{cor:DRpaths} and Lemma~\ref{l:lowerBound}, we now give the exact value of the game isolation number for certain paths.

\begin{corollary}
Let $n \ge 6$.  If $n \equiv 1 \pmod 5$, $n \equiv 2 \pmod 5$, or $n \equiv 3 \pmod 5$, then
\[
\io(P_n)=\io'(P_n)= \left\lfloor\frac{2n+2}{5}\right\rfloor\,.\]
\end{corollary}

\proof

When $n \equiv 1 \pmod 5$ or $n \equiv 3 \pmod 5$ the lower bound and the upper bound in Corollary~\ref{cor:DRpaths} agree.  Therefore, we get
$\io(P_n)=\io'(P_n)= \left\lfloor\frac{2n+2}{5}\right\rfloor$.

Finally, let $n \equiv 2 \pmod 5$, i.e. $n=5k+2$ for some $k \in \mathbb{N}$. By Corollary~\ref{cor:DRpaths}, it is enough to prove that $\io(P_n) \geq \left\lfloor\frac{2n+2}{5}\right\rfloor.$ For the purpose of contradiction assume that $\io(P_n) \leq \left\lfloor\frac{2n+2}{5}\right\rfloor -1 = 2k$. By Lemma~\ref{l:lowerBound}, we get $\io(P_n)=2k$. Now we use exactly the same strategy of Staller and the same notation as in the proof of Lemma~\ref{l:lowerBound}.
Recall that $p$ is the number of moves played in the D-game when Staller employs the above strategy and $p=\sum_{i=1}^k{|B_i|}$. Now, by our assumption $p=\io(P_n)=2k$. This implies that the last move was done by Staller. Using notation from Lemma~\ref{l:lowerBound} we get $p=\sum_{i=1}^k|B_i|=2k$ and $\sum_{i=1}^{k+1}|A_i| \leq 3k+1$. Thus $n \leq 5k+1$, which gives $2k=p \geq \frac{2n-2}{5}=\frac{2(5k+2)-2}{5}=2k+ \frac{2}{5},$ a contradiction.

\qed

\medskip

\section*{Acknowledgements}
B.B. and T.D. acknowledge the financial support of the Slovenian Research and Innovation Agency (research core funding No.\ P1-0297 and projects N1-0285, J1-3002, and J1-4008). All five
authors thank the Slovenian Research and Innovation Agency for financing our bilateral project between Slovenia and the
USA (title: Domination concepts in graphs, project No. BI-US/24-26-036).

\end{document}